\documentclass[a4paper,12pt]{article}
\usepackage{amssymb}

\newcommand{\R}{\mathbb{R}}

\newcommand{\N}{\mathbb{N}}

\newcommand{\cuad}{{\sqcap\kern-.68em\sqcup}}

\newcommand{\norm}[1]{\|#1\|}

\newtheorem{definition}{Definition}[section]
\newtheorem{theorem}{Theorem}[section]
\newtheorem{proposition}{Proposition}[section]

\newtheorem{lemma}{Lemma}[section]
\newtheorem{corollary}{Corollary}[section]
\newtheorem{remark}{Remark}[section]
\newcommand{\bremark}{\begin{remark} \em}
\newcommand{\eremark}{\end{remark} }

\begin{document}

\begin{center}{\bf  \large    Isolated singularities for  elliptic  equations with Hardy\\[2mm]

 operator and source nonlinearity  }\medskip
\bigskip\medskip

 {\small  Huyuan Chen\footnote{chenhuyuan@yeah.net} \quad and \quad  Feng Zhou\footnote{ fzhou@math.ecnu.edu.cn}
}
 \bigskip\medskip

{\small  $ ^1$Department of Mathematics, Jiangxi Normal University,\\
Nanchang, Jiangxi 330022, PR China\\[3mm]

$ ^2$Center for PDEs and Department   of Mathematics,
 East China Normal University,\\    Shanghai, 200241,  PR China
} \\[6mm]

 \begin{abstract}
 In this paper,  we concern the isolated singular solutions for semi-linear elliptic equations involving  Hardy-Leray potential
 \begin{equation}\label{0}
 -\Delta u+\frac{\mu}{|x|^2} u=u^p\quad {\rm in}\quad \Omega\setminus\{0\},\qquad u=0\quad{\rm on}\quad \partial\Omega.
 \end{equation}
 We classify the isolated singularities and obtain the existence and stability of positive solutions of (\ref{0}). Our results are based on
 the study of nonhomogeneous Hardy problem in a new distributional sense.
\end{abstract}

\end{center}
  \noindent {\small {\bf Key Words}:   Hardy potential, Isolated singularity,  Classification.}\vspace{1mm}

\noindent {\small {\bf MSC2010}:     35B44,     35J75. }

\vspace{2mm}

\setcounter{equation}{0}
\section{Introduction}

It is well known that the fundamental solutions of Laplacian operator play an essential role in the study of isolated singularities
of semilinear elliptic equations. Brezis and Lions in \cite{BL}   made use of the Schwartz theorem to show that any nonnegative solution of
$-\Delta u=f$ in $\Omega\setminus \{0\}$ must be a distributional solution of $-\Delta u=f+k\delta_0$ with $k\ge0$ and $\delta_0$ is the Dirac mass concentrated at the origin.
As a consequence, the solution would behavior as multiple of fundamental solutions. Later on, Lions in \cite{L} made use of this observation to
classify the isolated singularities of semilinear elliptic equation
\begin{equation}\label{eq 1.5}
-\Delta u=u^p\quad{\rm in}\quad \Omega\setminus\{0\},\qquad u=0\quad{\rm on}\quad \partial\Omega,
\end{equation}
by building the connection with distributional solution of
\begin{equation}\label{eq 1.5.1}
-\Delta u=u^p+k\delta_0\quad{\rm in}\quad \Omega.
\end{equation}
This work shows that when $p\ge \frac{N}{N-2}$, then $k=0$; when
$1<p<\frac{N}{N-2}$, there exists $k^*>0$ such that  (\ref{eq 1.5.1}) has two solutions for $k\in(0,k^*)$, which are classical solutions of (\ref{eq 1.5}), \ does one for $k=k^*$ and
no solution for $k>k^*$. Here $\frac{N}{N-2}$ is called as
Serrin's exponent. Motivated by the study of Lions' work,  Naito-Sato in \cite{NS} obtained two solutions for the model
$$-\Delta u+u=u^p+k_i\sum_{x_i}\delta_{x_i}\quad {\rm in}\quad \R^N,$$
by variational methods.

Note that when $p\ge\frac{N}{N-2}$ ,  the isolated singularity of (\ref{eq 1.5}) is invisible in the distributional sense  by Dirac mass as in   (\ref{eq 1.5.1}).
Using dynamic analysis, the classification of positive singularities of (\ref{eq 1.5})  has been done  by Aviles in \cite{A} for $p = \frac{N}{N-2}$ and $N\ge 3$, by Gidas and  Spruck in \cite{GS}
for $\frac{N}{N-2} < p < \frac{N+2}{N-2}$, by Caffarelli, Gidas and  Spruck in \cite{CGS} for
$p=\frac{N+2}{N-2}$. The books \cite{GT,V} give a survey on positive singularities of semilinear elliptic equations. Using this classification, solutions of equation (\ref{eq 1.5}) with many singular points were constructed in \cite{MP,P1}.
Recently,  in \cite{CZ1,CZ2}, we developed the Lions' method to classify isolated singularities of   Choquard equation in the subcritical case.

During the last years, there has been a renewed and increasing
interest in the study of linear and nonlinear elliptic equations involving Hardy operator,
 denoting
 $$\mathcal{L}_\mu=-\Delta+\mu|x|^{-2},$$
 motivated by great
applications and  important advances on the theory of nonlinear
partial differential equations, for instant, \cite{CL,DD,D,FF,JL,PT} and the references therein.
 The main tool to derive solutions is the variational methods due to the Hardy type inequality, see \cite{ACR,BM,BV1,FT,GP}.

When $\mu\ge\mu_0:=-\frac{(N-2)^2}{4}$, by Hardy inequality and Lax-Milgram theorem, it is known that
for $f\in L^2(\Omega)$,   non-homogeneous Hardy problem
 \begin{equation}\label{eq 1.1f}
 \displaystyle   \mathcal{L}_\mu u= f\qquad
   {\rm in}\quad  {\Omega}\setminus \{0\},\qquad u=0\quad{\rm on}\quad \partial\Omega
\end{equation}
has a unique solution  $u\in H^1_0(\Omega)$. Here and in the sequel, we always assume that $\Omega$
is a bounded, smooth domain containing the origin. A natural question raised:  what is the sharp condition of $f$ for the existence or nonexistence of (\ref{eq 1.1f})?
An attempt done by Dupaigne in \cite{D} is to consider  problem (\ref{eq 1.1f}) when $\mu\in[\mu_0,0)$ and $N\ge3$   in the distributional sense,
\begin{equation}\label{d1}
 \int_\Omega u\mathcal{L}_\mu \xi\, dx=\int_\Omega f\xi\, dx,\quad\forall\,\xi\in C^\infty_c(\Omega).
\end{equation}
 When $N\ge3$ and $\mu\in[\mu_0,0)$, the corresponding semi-linear problem has been studied by \cite{BP,FM1}.

While  a  critical defect  is that the singularity of the fundamental solution $\Phi_\mu,\Gamma_\mu$  of $\mathcal{L}_\mu$  is invisible from  (\ref{d1}) with $\Omega=\R^N$,  and   even $\Phi_\mu$  of $\mathcal{L}_\mu$ fails for (\ref{d1}) with $\Omega=\R^N$, when $\mu>0$,
where
 $$\Phi_\mu(x)=\arraycolsep=1pt\left\{
\begin{array}{lll}
|x|^{\tau_-(\mu)}\quad&{\rm if}\quad \mu>\mu_0\\[1mm]
|x|^{\tau_-(\mu)}(-\ln|x|) \quad&{\rm if}\quad \mu=\mu_0
\end{array}\right.
\quad{\rm and}\quad \Gamma_\mu(x)=|x|^{ \tau_+(\mu)},
$$
are two radially symmetric solutions of  problem
\begin{equation}\label{eq 1.02}
  \mathcal{L}_\mu u= 0\quad{\rm in}\quad  \R^N\setminus \{0\}.
 \end{equation}
Here
\begin{equation}\label{1.1}
 \tau_\pm(\mu)=- (N-2)/2\,\pm\sqrt{\mu-\mu_0}
\end{equation}
 are two roots of $\mu-\tau(\tau+N-2)=0$.

To overcome this defect arising from the Hardy potential,   a new distributional identity has been proposed in \cite{CQZ} recently.
For $\mu\ge\mu_0$, $\Phi_\mu $ could be seen as a $d\mu$-distributional solution of
\begin{equation}\label{eq003}
\mathcal{L}_\mu   u =b_{\mu}\delta_0 \quad   {\rm in}\quad \mathcal D'(\R^N),
\end{equation}
  in the distribution sense
\begin{equation}\label{eq0030}
\int_{\R^N}\Phi_\mu \mathcal{L}_\mu^*\xi\, d\mu =b_\mu\xi(0),\quad\forall\, \xi\in C_c^\infty(\R^N),
\end{equation}
with the measure $d\mu(x) =\Gamma_\mu(x) dx$ and
\begin{equation}\label{L}
\mathcal{L}^*_\mu=-\Delta -2\frac{\tau_+(\mu) }{|x|^2}\,x\cdot\nabla
\end{equation}
and the normalized constant
$$
b_{\mu}: =\left\{\arraycolsep=1pt
\begin{array}{lll}
 2\sqrt{\mu-\mu_0}\,|\mathcal{S}^{N-1}|\quad
   &{\rm if}\quad \mu>\mu_0,\\[1.5mm]
 \phantom{   }
|\mathcal{S}^{N-1}| \quad  &{\rm  if}\quad \mu=\mu_0,
 \end{array}
 \right.
$$
where $\mathcal{S}^{N-1}$ is the unit sphere in $\R^N$ and $|\mathcal{S}^{N-1}|$ is the area of the unit sphere.

For problem (\ref{eq 1.1f}),  the $d\mu$-distributional sense provides a complete understanding of  the existence, non-existence and the singularities of isolated singular
solutions.

\begin{theorem}\label{teo 2.1}\cite[Theorem 1.3]{CQZ}
Let    $f$ be a function in $C^\gamma_{loc}(\overline{\Omega}\setminus \{0\})$ for some $\gamma\in(0,1)$.

$(i)$ Assume that
\begin{equation}\label{f1}
\int_{\Omega} |f|\,   d\mu <+\infty,
\end{equation}
  then    problem (\ref{eq 1.1f}), subject to $\lim_{x\to0}u(x)\Phi_\mu^{-1}(x)=k$ with $k\in\R$,
has a  unique solution $u_k$, which  satisfies the distributional identity
 \begin{equation}\label{1.2f}
 \int_{\Omega}u_k  \mathcal{L}_\mu^*(\xi)\, d\mu  = \int_{\Omega} f  \xi\, d\mu +b_\mu k\xi(0),\quad\forall\, \xi\in   C^{1.1}_0(\Omega).
\end{equation}

$(ii)$ Assume that $f$ verifies (\ref{f1})  and $u$ is a nonnegative solution of (\ref{eq 1.1f}), then  $u$ satisfies (\ref{1.2f}) for some $k\ge0$
and verifies that $\lim_{x\to0}u(x)\Phi_\mu^{-1}(x)=k$.\smallskip

$(iii)$ Assume that $f\ge0$ and
\begin{equation}\label{f2}
 \lim_{r\to0^+} \int_{\Omega\setminus B_r(0)} f\, d \mu  =+\infty,
\end{equation}
then problem (\ref{eq 1.1f}) has no nonnegative solutions.
\end{theorem}

Our concern of this article is  to analyze  the isolated singular solution of semi-linear  problem
\begin{equation}\label{eq 1.1}
    \mathcal{L}_\mu u= u^p\quad
   {\rm in}\quad  \Omega\setminus \{0\},\qquad u=0\quad{\rm on}\quad \partial\Omega,
 \end{equation}
where    $p>1$ and  $\Omega$ is a bounded $C^2$ domain containing the origin  in $\mathbb{ R}^N$.

 Our first result on the classification of isolated singularities of (\ref{eq 1.1}), based on Theorem \ref{teo 2.1}, states as following.
\begin{theorem}\label{teo 0}
Assume that $p>1$, $d\mu$ and $\mathcal{L}_\mu^*$ are given in (\ref{L}), $\tau_\pm(\mu)$ are given in (\ref{1.1}) and $u$ is a nonnegative classical solution of (\ref{eq 1.1}) in $\Omega \setminus \{0\}$. Let
$$p^*_\mu=1+\frac{2}{ -\tau_-(\mu)}.$$

     Then  $u\in L^p(\Omega,d\mu)$ and
there exists $k\ge0$   such that  $u$ is a $d\mu$-distributional solution of
 \begin{equation}\label{eq 1.2}
    \displaystyle    \mathcal{L}_\mu    u-u^p =k\delta_0\quad
 {\rm in}\quad \Omega, \qquad u=0\quad{\rm on}\quad \partial\Omega
\end{equation}
that is,
 \begin{equation}\label{e 1.1}
     \int_\Omega \left(u \mathcal{L}^*_\mu \xi-u^p\xi\right)\, d\mu=k\xi(0),\quad \forall \xi\in C^\infty_c(\Omega).
\end{equation}
Furthermore,
$(i)$  when $p\ge p^*_\mu$, we have that $k=0$;

$(ii)$ Assume that $p\in(1,p^*_\mu)$,  if $k>0$, then $u$    satisfies that
\begin{equation}\label{1.2}
 \lim_{x\to0} u(x)\Phi_\mu^{-1}(x)=b_\mu  k;
\end{equation}
and if $k=0$, then  $u$    satisfies that
\begin{equation}\label{1.2-1}
 \limsup_{x\to0} u(x)\Gamma_\mu^{-1}(x)<+\infty.
\end{equation}

\end{theorem}

Note that the isolated singularities of nonnegative solutions of (\ref{eq 1.1}) verifying (\ref{1.2}) with $k>0$ could be
seen by Dirac mass in the subcritical case $1< p<p^*_\mu$. While the singularities is invisible in the super critical case $p\ge p^*_\mu$, which also shows that
the singularities is weaker that $\Phi_\mu$. We note also that if $\mu =0$, then $p^*_\mu= \frac{N}{N-2}$ and it recovers the early result of Lions in \cite{L}. It is open but challenging  to clarify the singularities in the super critical case   $p\ge p^*_\mu$.

Concerning the existence of singular solutions of (\ref{eq 1.1}), we can prove the following:
\begin{theorem}\label{teo 1}
Assume that $p\in(1,p^*_\mu)$,  then there exists $k^*>0$
  such that

  $(i)$ for $ k\in(0, k^*)$, problem (\ref{eq 1.1}) subjecting to (\ref{1.2})  admits a minimal nonnegative solution $u_{k}$
and  a Mountain Pass type solution $w_k>u_k$, both are  $d\mu$-distributional solutions of (\ref{eq 1.2});

$(ii)$ for $k=k^*$, problem (\ref{eq 1.2})  admits a $d\mu$-distributional solution $u_{k^*}$. In the particular case that  $\Omega=B_1(0)$,    problem (\ref{eq 1.1}) subjecting to (\ref{1.2}) with $k=k^*$ admits a unique solution  $u_{k^*}$, which is a very weak solution of (\ref{eq 1.2});

$(iii)$
for $k> k^*$,  problem (\ref{eq 1.1}) subjecting to (\ref{1.2})  admits no solution.

\end{theorem}

 When $\mu=0$, Theorem \ref{teo 1} could be seen in \cite{L}. Given some $k$, the minimal solution $u_k$ of (\ref{eq 1.2}) is obtained by the following iterating procedure:
$$v_0= k \mathbb{G}_\Omega [\delta_0], \qquad v_n  = \mathbb{G}_\Omega [v_{n-1}^p]+ k \mathbb{G}_\Omega [\delta_0],$$
where $\mathbb{G}_\Omega$ is the Green's operator defined by the Green kernel $G(x,y)$ of $  -\Delta $ in $\Omega\times\Omega$ under the zero Dirichlet
boundary condition. To control the sequence $\{v_n\}_n$,  a barrier function is constructed by the estimate that
 \begin{equation}\label{1.5q}
 \mathbb{G}_\Omega[\mathbb{G}_\Omega^p[\delta_0]]\le c_2\mathbb{G}_\Omega[\delta_0]\quad {\rm in}\quad \Omega\setminus\{0\}.
 \end{equation}
The  $d\mu$-distributional solution could be improved into the classical solution of (\ref{eq 1.1}) when $\mu=0$. While for $\mu\ge \mu_0$ but $\mu\not=0$, it
is difficult to improve the very weak solution of (\ref{eq 1.2}) to be a classical solution of (\ref{eq 1.1}). To overcome this difficulty,
our idea is to derive the minimal solution of (\ref{eq 1.1}) subjecting to (\ref{1.2}) with $k\in(0,k^*)$ directly.  The optimal value of $k$ for existence is given by
$$k^*=\sup\{k>0:\ (\ref{eq 1.1})\ {\rm subjecting\ to\  (\ref{1.2}) \ with\ such}\ k \ {\rm  has\ minimal\ solution} \},$$
which is equal to or less than the optimal value  for existence of weak solutions to  (\ref{eq 1.2}). This gives rise to a difficulty for the uniqueness when $k=k^*$.
In the particular case of $\Omega=B_1(0)$, we make use of the properties of the radial symmetry and the monotonicity, we can get the uniqueness and improve the regularity for $k=k^*$.

The rest of this paper is organized as follows. In Section 2, we  show  the Comparison Principle, which  is important for the existence of isolated solutions
  of (\ref{eq 1.1}) under the constraint of (\ref{1.2}) with $k\in(0,k^*)$ and we do the classification of singularities of positive solutions for (\ref{eq 1.1}).
 Section 3 is devoted to isolated singular   solutions of (\ref{eq 1.1}) in the subcritical case.

\vskip0.5cm

\setcounter{equation}{0}
\section{Proof of Theorem \ref{teo 0}}

In this section, we concentrate on the classification of  isolated singular solutions to (\ref{eq 1.1}). In what follows,  we denote by $c_i$ a generic  positive constant.
We first introduce some basic tools in the classical sense. One basic tool is the Comparison Principle.

\begin{lemma}\label{lm cp}
Let $O$ be a bounded open set in $\R^N$, $L: O\times [0,+\infty)\to[0,+\infty)$ be a continuous function satisfying that for any $x\in  O$,
$$L(x,s_1)\ge L(x,s_2)\quad {\rm if}\quad s_1\ge s_2,$$ then $\mathcal{L}_\mu+L$ with $\mu\ge \mu_0$ verifies the Comparison Principle,
that is, if
$$u,\,v\in C^{1,1}(O)\cap C(\bar O)$$ verify that
$$\mathcal{L}_\mu u+ L(x,u)\ge \mathcal{L}_\mu v+ L(x,v) \quad {\rm in}\quad  O
\qquad{\rm and}\qquad   u\ge  v\quad {\rm on}\quad \partial O,$$
then
$$u\ge v\quad{\rm in}\quad  O.$$

\end{lemma}
{\bf Proof.} Let $w=u-v$ and then $w\ge 0$ on  $\partial O$.
Let $w_-=\min\{w,0\}$ and our purpose is to prove that $w_-\equiv 0$.
If $ O_-:=\{x\in O:\, w(x)<0\}$ is not empty, then it is a bounded
$C^{1,1}$ domain in $ O$ and $w_-=0$ on $\partial O$ by the assumption that  $u\ge v$ on $\partial O$.  From Hardy inequality with $\mu\ge \mu_0$,  there holds,
\begin{eqnarray*}
 0  &=& \int_{ O_-}(-\Delta  w_- +\frac{\mu}{|x|^{2 }}  w_- ) w_- dx+\int_{ O_-}[L(x,u)-L(x,v)](u-v)_-\,dx \\ &\ge&  \int_{ O_-} \left(|\nabla w_-|^2 +\frac{\mu}{|x|^{2 }}  w_-^2\right) dx
 \ge  c_1 \int_{ O_-}w_-^2 dx,
\end{eqnarray*}
then  $w_-=0$ in a.e. $ O_-$, which is impossible with the definition of $ O_-$.\hfill$\Box$\smallskip

As a consequence, we have the following lemma which plays an important role in the obtention of the uniqueness for classical solution.

\begin{lemma}\label{cr hp}
Assume that $\Omega$ is a bounded $C^2$ domain and $u_i$ with $i=1,2$ are classical solutions of
\begin{equation}\label{eq0 2.1}
 \arraycolsep=1pt\left\{
\begin{array}{lll}
 \displaystyle \mathcal{L}_\mu u = f_i\qquad
   {\rm in}\quad  {\Omega}\setminus \{0\},\\[1.5mm]
 \phantom{ L_\mu     }
 \displaystyle  u= 0\qquad  {\rm   on}\quad \partial{\Omega},\\[1.5mm]
 \phantom{   }
  \displaystyle \lim_{x\to0}u(x)\Phi_\mu^{-1}(x)=k
 \end{array}\right.
\end{equation}
for some $k\in \R$. If $f_1\le f_2$ in $\Omega\setminus \{0\}$, then
$$u_1\le u_2\quad{\rm in}\quad \Omega\setminus \{0\}.$$

\end{lemma}
{\bf Proof.}  Let $u=u_1-u_2$ satisfy that
$$\mathcal{L}_\mu u \le 0,$$
 then for any $\epsilon>0$, there exists $r_\epsilon>0$ converging to zero as $\epsilon\to0$ such that
 $$u\le \epsilon \Phi_\mu\quad{\rm in}\quad \overline{B_{r_\epsilon}(0)}\setminus\{0\}.$$
We see that
$$u=0<\epsilon \Phi_\mu \quad{\rm on}\quad \partial\Omega,$$
then by Lemma \ref{lm cp}, we have that
$$u\le \epsilon \Phi_\mu\quad{\rm in}\quad \Omega\setminus\{0\}. $$
By the arbitrary of $\epsilon$, we have that $u\le 0$  in $\Omega\setminus\{0\}$, which ends the proof.
\hfill$\Box$\medskip

\begin{lemma}\label{lm 2.1-singular}
Assume that  $\mu>\mu_0$ and $f\in C^1(\Omega\setminus \{0\})$ verifies that
\begin{equation}\label{2.1}
0\le f(x)\le c_2 |x|^{\tau-2},
\end{equation}
where   $\tau>\tau_-(\mu)$.
Let $u_f$ be the solution of
\begin{equation}\label{2.2}
 \arraycolsep=1pt\left\{
\begin{array}{lll}
 \displaystyle  \mathcal{L}_\mu u= f\qquad
   {\rm in}\quad \ {\Omega}\setminus \{0\},\\[1mm]
 \phantom{  L_\mu \, }
 \displaystyle  u= 0\qquad  {\rm   on}\quad\ \partial{\Omega},\\[1mm]
 \phantom{   }
  \displaystyle \lim_{x\to0}u(x)\Phi_\mu^{-1}(x)=0.
 \end{array}\right.
\end{equation}
Then we have that if $\tau_-(\mu)<\tau<\tau_+(\mu)$,
\begin{equation}\label{2.3}
0\le u_f(x)\le c_3|x|^{\tau}\qquad    {\rm in}\quad \ {\Omega}\setminus \{0\};
\end{equation}
 if $ \tau=\tau_+(\mu)$,
\begin{equation}\label{2.3-1}
0\le u_f(x)\le c_4|x|^{\tau}(1+(-\ln |x|)_+)\qquad    {\rm in}\quad \ {\Omega}\setminus \{0\};
\end{equation}
and if $\tau>\tau_+(\mu)$,
\begin{equation}\label{2.3-2}
0\le u_f(x)\le c_5|x|^{\tau_+(\mu)} \qquad    {\rm in}\quad \ {\Omega}\setminus \{0\}.
\end{equation}

\end{lemma}
{\bf Proof.} We only have to construct suitable upper bound for $u_f$.
For $\tau_-(\mu)<\tau< \tau_+(\mu)$, the upper bound is $t|x|^{\tau}$ for some suitable $t$.
In fact, $$\mathcal{L}_\mu |x|^{\tau}=c_\tau |x|^{\tau-2},$$
where $c_\tau>0$. So for some $t>0$, there holds
$$\mathcal{L}_\mu (t|x|^{\tau}) \ge f(x),\qquad    \forall\, x\in {\Omega}\setminus \{0\}.$$
Then (\ref{2.3}) follows by Lemma \ref{cr hp}.

For  $\tau= \tau_+(\mu)$, we have that
$$\mathcal{L}_\mu (|x|^{\tau_+(\mu)}(-\ln |x|))=(2\tau_+(\mu)+N-2)|x|^{\tau_+(\mu)}. $$
The upper bound could be constructed by adjusting the coefficients of
$$ s|x|^{\tau_+(\mu)}+t|x|^{\tau_+(\mu)}(-\ln |x|).$$

For  $\tau>\tau_+(\mu)$, we have that
$$\mathcal{L}_\mu (|x|^\tau)= (\mu-\tau(\tau+N-2))|x|^{\tau-2}, $$
where $\mu-\tau(\tau+N-2)<0$ for $\tau>\tau_+(\mu)$.
The upper bound could be constructed by adjusting the coefficients $s,t$ of
$$ s|x|^{\tau_+(\mu)}-t|x|^{\tau}.$$

\medskip

\noindent{\bf Proof of Theorem \ref{teo 0}.} Let $v$ be a nonnegative nontrivial solution of (\ref{eq 1.1}), then $v^p$ is a nonnegative nontrivial function.
Taking $f=v^p$ in (\ref{eq 1.1f}),  we have that $v^p\in L^1(\Omega,\,d\mu)$, otherwise, from Theorem \ref{teo 2.1} part $(iii)$, we have that
there is no solution for problem
$$\arraycolsep=1pt\left\{
\begin{array}{lll}
 \displaystyle   \mathcal{L}_\mu u= v^p\qquad
   &{\rm in}\quad  {\Omega}\setminus \{0\},\\[1.5mm]
 \phantom{   L_\mu   }
 \displaystyle  u= 0\qquad  &{\rm   on}\quad \partial{\Omega},
 \end{array}\right.
 $$
which contradicts our assumption.

From Theorem \ref{teo 2.1} part $(i)$ and $(ii)$, we know that $v$ is a weak solution of (\ref{eq 1.2}) for some $k\ge0$ and
\begin{equation}\label{2.4}
 \lim_{x\to0}v(x)\Phi_\mu^{-1}(x)=b_\mu k.
\end{equation}

For $p\ge  p^*_\mu$, we show that $k=0$. If not, assume that $k>0$ and there exists $r_0>0$ such that
$$v(x)\ge \frac k2 \Phi_\mu(x),\quad\forall\, x\in  B_{r_0}(0)\setminus\{0\},$$
which implies that
$$v^p(x)\Gamma_\mu(x)\ge \frac k2 |x|^{p\tau_-(\mu)+\tau_+(\mu)},\quad\forall\, x\in  B_{r_0}(0)\setminus\{0\},$$
where
$$p\,\tau_-(\mu)+\tau_+(\mu)\le  p^*_\mu\,\tau_-(\mu)+2-N-\tau_-(\mu)= -N, $$
so $v^p\not\in L^1(\Omega,\,d\mu)$, that contradicts  $v^p\in L^1(\Omega,\,d\mu)$.

For $1<p<  p^*_\mu$ and $k>0$, (\ref{1.2}) follows by  (\ref{2.4}).

For  $1<p<  p^*_\mu$ and $k=0$, (\ref{2.4}) implies that
$$\lim_{x\to0}v(x)\Phi_\mu^{-1}(x)=0.$$
We have that for some $d_0>0$,
$$v^p(x)\le d_0 \Phi_\mu^p(x),\quad\forall\, x\in \Omega\setminus\{0\}, $$
and then
$$v^p(x)\le d_0|x|^{\tau_1-2},\quad\forall\, x\in \Omega\setminus\{0\} $$
where
$\tau_1=\tau_-(\mu) p+2-\epsilon$
and $\epsilon>0$ is such that $\tau_-(\mu) p+2-\epsilon>\tau_-(\mu)$,
$\tau_1\not=\tau_+(\mu).$
From Lemma \ref{lm 2.1-singular}, we have that
$$v(x)\le d_1|x|^{\min\{\tau_+(\mu),\, \tau_1\} },\quad\forall\, x\in \Omega\setminus\{0\}. $$
If $\tau_1>\tau_+(\mu)$, we are done. Especially, for $\mu=\mu_0$, we have
$$\tau_1>\tau_-(\mu)=\tau_+(\mu)$$
and then (\ref{1.2-1}) holds.

If not, we only consider the case $\mu>\mu_0$, let $\tau_2=\tau_1 p+2$ and by adjusting $\epsilon$, we have that
$$\tau_2\not=\tau_+(\mu)$$
and
$$v(x)\le d_2|x|^{\min\{\tau_+(\mu),\, \tau_2\} },\quad\forall\, x\in B_{r_0}(0)\setminus\{0\}. $$

Let
$$\tau_j=\tau_{j-1}p+2,\quad j=2,3\cdots,$$
which is  an increasing sequence s and
$$\lim_{j\to+\infty}\tau_j=+\infty.$$
So there exists $j_0$ such that
$$\tau_{j_0}\ge \tau_+(\mu)\quad{\rm and}\quad \tau_{j_0-1}<\tau_+(\mu).$$
Adjusting $\epsilon$, we can improve that
$$\tau_{j_0}> \tau_+(\mu)\quad{\rm and}\quad \tau_{j_0-1}<\tau_+(\mu).$$
We observe that
for $j\le j_0$,
$$v(x)\le d_{j-1} |x|^{\tau_{j-1}},  $$
then
$$v(x)\le d_{j} |x|^{\tau_{j}}.  $$
By Lemma \ref{lm 2.1-singular} iteratively until $j=j_0$, we can obtain that  (\ref{1.2-1}) holds.
\hfill$\Box$

\setcounter{equation}{0}
\section{Existence and Stability }

In this section, we search for the singular solutions of problem (\ref{eq 1.1}) provided the asymptotic behavior at the origin as $\lim_{|x|\to0} u(x)\Phi_\mu^{-1}(x)=b_\mu k$.

\subsection{Existence of Minimal Solution}

\begin{proposition}\label{pr m1}
Assume that $p\in(1,p^*_\mu)$,  then there exists $k^*>0$
  such that

 $(i)$ for $ k\in(0, k^*)$, problem (\ref{eq 1.1}) subjecting to (\ref{1.2})  admits a minimal nonnegative solution $u_{k}$; and $u_k$
 is a very weak solution of (\ref{eq 1.2}).

$(ii)$
for $k> k^*$,  problem (\ref{eq 1.1}) subjecting to (\ref{1.2})  admits no solution.

\end{proposition}



\noindent {\bf Proof.} For $k>0$, let $v_0$ be the solution of
\begin{equation}\label{2.02-1}
 \arraycolsep=1pt\left\{
\begin{array}{lll}
 \displaystyle  \mathcal{L}_\mu u= 0\qquad
   {\rm in}\quad \ {\Omega}\setminus \{0\},\\[1mm]
 \phantom{  L_\mu \, }
 \displaystyle  u= 0\qquad  {\rm   on}\quad\  \partial{\Omega},\\[1mm]
 \phantom{   }
  \displaystyle \lim_{x\to0}u(x)\Phi_\mu^{-1}(x)=b_\mu k
 \end{array}\right.
\end{equation}
and from Theorem \ref{teo 0} part (iii),
we may define the iterating sequence $ v_n$, the solution of
\begin{equation}\label{2.02-2}
 \arraycolsep=1pt\left\{
\begin{array}{lll}
 \displaystyle  \mathcal{L}_\mu u= v_{n-1}^p\quad
   {\rm in}\quad  {\Omega}\setminus \{0\},\\[1mm]
 \phantom{  L_\mu \, }
 \displaystyle  u= 0\qquad \ {\rm   on}\quad \partial{\Omega},\\[1mm]
 \phantom{   }
  \displaystyle \lim_{x\to0}u(x)\Phi_\mu^{-1}(x)=b_\mu k.
 \end{array}\right.
\end{equation}
By Lemma \ref{cr hp}, a standard iteration argument shows that
$\{v_n\}_n$ is an increasing sequence of functions in $\Omega \setminus \{0\}$.

Let $w_0$ be the solution of (\ref{2.02-1}) with $k=1$ and $w_1$ be the solution of
\begin{equation}\label{2.02-3}
 \arraycolsep=1pt\left\{
\begin{array}{lll}
 \displaystyle  \mathcal{L}_\mu u= w_0^p\qquad
   {\rm in}\quad  {\Omega}\setminus \{0\},\\[1mm]
 \phantom{  L_\mu \, }
 \displaystyle  u= 0\qquad\ \,  {\rm   on}\quad \partial{\Omega},\\[1mm]
 \phantom{   }
  \displaystyle \lim_{x\to0}u(x)\Phi_\mu^{-1}(x)=0.
 \end{array}\right.
\end{equation}
We claim that
\begin{equation}\label{3.1}
w_1\le c_6w_0\quad    {\rm in}\quad  {\Omega}\setminus \{0\}.
\end{equation}
In fact, we see that $w_0,\,w_1$ are strictly positive and continuous in $\Omega\setminus\{0\}$,
$$\lim_{x\to0}w_0(x)\Phi_\mu^{-1}(x)=b_\mu,\quad \quad \lim_{x\to0}w_1(x)\Phi_\mu^{-1}(x)=0$$
and near the boundary,
$$\frac1{c_7}\le w_0(x)\rho^{-1}(x),\, w_1(x)\rho^{-1}(x)\le c_7,$$
thus,   (\ref{3.1}) holds.

Now we construct an upper bound for the sequence $\{v_n\}_n$ for suitable $k$. Let $ w_t$ be the solution of
\begin{equation}\label{2.02-4}
 \arraycolsep=1pt\left\{
\begin{array}{lll}
 \displaystyle  \mathcal{L}_\mu u= tk^pw_0^p\quad\quad
   {\rm in}\quad  {\Omega}\setminus \{0\},\\[1mm]
 \phantom{  L_\mu \, }
 \displaystyle  u= 0\qquad \ \ \quad {\rm   on}\quad \partial{\Omega},\\[1mm]
 \phantom{   }
  \displaystyle \lim_{x\to0}u(x)\Phi_\mu^{-1}(x)=b_\mu k
 \end{array}\right.
\end{equation}
and by (\ref{3.1}), we have that
$$w_t=t k^p w_1+kw_0\le (c_6tk^p+k) w_0,$$
so $w_t$ verifies that
\begin{equation}\label{2.02-5}
\mathcal{L}_\mu w_t\ge w_t^p,
\end{equation}
if
$$w_t^p\le (c_6tk^p+k)^p w_0^p\le tk^pw_0^p,$$
which holds for $t$ satisfying
\begin{equation}\label{2.5}
t\ge (c_6tk^{p-1}+1)^p.
\end{equation}
Note that the convex function $f_{k}(t) =  (c_6t k^{p-1} + 1)^p$  can intersect the line $g(t) = t$,  if
\begin{equation}\label{4.2.5}
    c_6 k^{p-1}\le \frac1p\left(\frac{p-1}p\right)^{p-1}.
\end{equation}
Let $k_p=\left(\frac1{c_6 p}\right)^{\frac1{p-1}}\frac{p-1}p$, then if $k\le k_p$, it always holds that $f_{k}(t_p)\le t_p$ for $ t_p=\left(\frac p{p-1}\right)^{p}.$
 Hence,  by the definition of $w_{t_p}$, we have that  $w_{t_p}>v_0$ and from Lemma \ref{cr hp}, it implies that
$$v_1 \le w_{t_p}.$$
Inductively, we obtain
\begin{equation}\label{2.10a}
v_n\le w_{t_p}\quad {\rm for\ all}\ n\in\N.
\end{equation}
 Therefore, the sequence $\{v_n\}_n$ converges. Let $u_{k}:=\lim_{n\to\infty} v_n$, then for any compact set $K$ in $\Omega\setminus\{0\}$,
 and then $u_k$ verifies the equation $$\mathcal{L}_\mu u= u^p\quad{\rm in}\quad K$$
and $v_0\le u_k\le w_{t_p}$, so  $u_{k}$ is a classical solution of (\ref{eq 1.1})   verifying (\ref{1.2})
and $u_k$ is a very weak solution of (\ref{eq 1.2}) with such $k$.

We claim that $u_{k}$ is the minimal solution of (\ref{eq 1.1}) verifying (\ref{1.2}), that is, for any nonnegative solution $u$ of (\ref{eq 1.1}) verifying (\ref{1.2}), we always have $u_{k}\leq u$. Indeed, from Lemma \ref{cr hp},  there holds
\[
 u  \ge v_0\quad{\rm and}\quad u^p\ge v_0^p,
\]
 then $ u \ge v_{1}.$
We may show inductively that
$
u\ge v_n
$
for all $n\in\N$.  The claim follows.

 Similarly,  if problem (\ref{eq 1.1}) subjecting to (\ref{1.2}) with $k=k_1>0$ has a nonnegative solution $u$, then (\ref{eq 1.1}) admits a minimal solution $u_{k}$ verifying (\ref{1.2}) with any $ k\in(0, k_1]$. As a result, the mapping $ k\mapsto u_{k}$ is increasing.
So we may define
$$k^*=\sup\{k>0:\ (\ref{eq 1.1})\ {\rm subjecting\ to\  (\ref{1.2}) \ with\ such}\ k \ {\rm  has\ minimal\ solution} \}.$$
We claim that $k^*<+\infty$. Observe that problem (\ref{eq 1.1}) subjecting to (\ref{1.2}) with $k>k^*$ has no solution and
$$k^*\ge k_p.$$

Now we prove $k^*<+\infty$.
 Let $(\lambda_1,\varphi_1)$ be the first eigenvalue and the corresponding nonnegative eigenfunction of $\mathcal{L}_\mu $ in $H^1_0(\Omega)$, see reference \cite{BM}. Taking
$\varphi_1=\Gamma_\mu \varphi_1^*$,
we have that
$$\Gamma_\mu  \mathcal{L}_\mu^*\varphi_1^*=  \mathcal{L}_\mu\varphi_1= \lambda_1\Gamma_\mu \varphi_1^*,\quad {\rm in}\quad \Omega\setminus\{0\}. $$
Assume that for $k>0$, (\ref{eq 1.1}) has a positive solution $u_k$   verifying (\ref{1.2}) with such $k$.
From Theorem \ref{teo 0}, $u_k$ is a very weak solution of (\ref{eq 1.2}). We next show that $k$ must be bounded. Let $\{\varphi_\epsilon^*\}_\epsilon$ be an increasing  sequence of nonnegative functions in $C^2_0(\Omega)$  such that
\begin{equation}\label{3.4}
 \varphi_\epsilon^*=\varphi_1^*\quad{\rm in}\quad \Omega\setminus B_\epsilon(0)\quad{\rm and}\quad \varphi_\epsilon^*\le \varphi_1^*\quad{\rm in}\quad   B_\epsilon(0),
\end{equation}
where $\epsilon>0$ will be determined latter.

By H\"{o}lder inequality, it implies that
\begin{eqnarray}
\lambda_1 \left(\int_{\Omega } u_{k}^p \varphi_\epsilon^*\,d\mu\right)^{\frac1p} \left(\int_{\Omega }  \varphi_\epsilon^*\,d\mu \right)^\frac{p-1}p  &\ge& \lambda_1 \int_{\Omega \setminus B_\epsilon(0)} u_{k} \varphi_\epsilon^*\,d\mu
\nonumber \\&=& \int_{\Omega} u_{k}\mathcal{L}_\mu^*(\varphi_\epsilon^*)\,d\mu-\int_{B_\epsilon(0)} u_{k} \mathcal{L}_\mu^*(\varphi_\epsilon^*)\,d\mu
\nonumber \\& =  &\int_{\Omega}   u_{k}^p\varphi_\epsilon^* \,d\mu+k\varphi_\epsilon^*(0)-\int_{B_\epsilon(0)} u_{k} \mathcal{L}_\mu^*(\varphi_\epsilon^*)\,d\mu
\nonumber\\& \ge& \int_{\Omega}   u_{k}^p\varphi_\epsilon^* \,d\mu-\int_{B_\epsilon(0)} u_{k} \mathcal{L}_\mu^*(\varphi_\epsilon^*)\,d\mu.\label{3.3}
\end{eqnarray}
From the fact that $\lim_{|x|\to 0} u_k(x)\Phi_\mu^{-1}(x)=b_\mu k$, there exists $\epsilon_0\in(0,1)$ small such that
 $u_k\le 2k\Phi_\mu$ in $B_{\epsilon_0}(0)\setminus\{0\}$. Then
\begin{eqnarray*}
|\int_{B_\epsilon(0)} u_{k} \mathcal{L}_\mu^*(\varphi_\epsilon^*)\,d\mu| &\le & c_8\int_{B_\epsilon(0)}\Phi_\mu\Gamma_\mu  |\mathcal{L}_\mu^*(\varphi_\epsilon^*)| dx
\\&\le & \arraycolsep=1pt\left\{
\begin{array}{lll}
c_8 \epsilon\quad    &{\rm for}\quad \mu>\mu_0,\\[1mm]
 \phantom{  }
c_8 \epsilon(-\ln \epsilon)\quad \ \ \quad &{\rm for}\quad \mu=\mu_0,
 \end{array}\right.
\end{eqnarray*}
which implies that
$$\lim_{\epsilon\to0}\int_{B_\epsilon(0)} u_{k} \mathcal{L}_\mu^*(\varphi_\epsilon^*)\,d\mu=0.$$
Since the mapping $\epsilon\to \varphi_\epsilon^*$ is decreasing, then  the mapping $\epsilon\to\int_{\Omega}   u_{k}^p\varphi_\epsilon^* \,d\mu$ is decreasing  and for any $ \upsilon\in(0,\frac12)$, there exists $\epsilon>0$ small such that
$$\int_{B_\epsilon(0)} u_{k} \mathcal{L}_\mu^*(\varphi_\epsilon^*)\,d\mu\le \upsilon \int_{\Omega}   u_{k}^p\varphi_\epsilon^* \,d\mu$$
Therefore, fixed $\upsilon=\frac14$, choosing $\epsilon>0$ suitable, we have that
$$\int_{\Omega}   u_{k}^p\varphi_1^* \,d\mu \le \left(\frac{4\lambda_1}{3}\right)^{\frac p{p-1}}\int_{\Omega}  \varphi_1^*\,d\mu.  $$
We observe that
$$u_k\ge kw_0 \quad{\rm in}\quad \Omega,$$
then
$$k^p \int_{\Omega}   w_0^p\varphi_\epsilon^* \,d\mu \le \left(\frac{4\lambda_1}{3}\right)^{\frac p{p-1}}\int_{\Omega}  \varphi_\epsilon^*\,d\mu, $$
which implies that
$$
 k\le \left(\frac{4\lambda_1}{3}\right)^{\frac 1{p-1}}\left(\frac{\int_{\Omega}  \varphi_\epsilon^*\,d\mu}{\int_{\Omega}   w_0^p\varphi_\epsilon^* \,d\mu}\right)^{\frac1p}.
$$
Therefore, $k^*\le \left(\frac{4\lambda_1}{3}\right)^{\frac1{p-1}}\left(\frac{\int_{\Omega}  \varphi_\epsilon^*\,d\mu}{\int_{\Omega}   w_0^p\varphi_\epsilon^* \,d\mu}\right)^{\frac1p}<+\infty$.
\hfill$\Box$

\subsection{Stability of Minimal Solution}

In this subsection, we discuss the stability of  minimal solution for problem (\ref{eq 1.1}).

\begin{definition}\label{def 1}
A solution (or weak solution) $u$ of (\ref{eq 1.1}) is stable (resp. semi-stable) if
$$
\norm{\xi}_\mu^2:=\int_\Omega(|\nabla \xi|^2+\mu|x|^{-2} \xi^2) dx > p\int_{\Omega}  u^{p-1}\xi^2 \,dx,\quad ({\rm resp.}\ \ge)\quad \forall \xi\in H^1_0(\Omega)\setminus\{0\}.
$$
We also note that $\norm{ \cdot}_\mu$ is a norm of $H^1_0(\Omega)$ induced  by the inner product for $\mu \ge \mu_0$,
$$\langle u,v\rangle_\mu=\int_\Omega( \nabla u\nabla v + \mu |x|^{-2} uv )\,dx.$$

\end{definition}

\begin{lemma}\label{lm stability}

For $k\in (0,k^*)$, let $u_{k}$ be the minimal nonnegative solution of (\ref{eq 1.1}) subjecting to (\ref{1.2}).
Then $u_{k}$ is stable.
Moreover,   for any $\xi\in H^1_0(\Omega)\setminus\{0\},$ we have that
\begin{equation}\label{e.3.1}
\norm{\xi}_\mu^2- p\int_{\Omega} u_{k}^{p-1}\xi^2 \,dx\ge  \left( 1-(\frac{2k}{k+k^*})^{p-1}\right) \norm{\xi}_\mu^2.
\end{equation}

\end{lemma}
{\bf Proof.} We first  prove the  stability of $u_k$ when $k>0$ small.
When $0<k<k_p$,  the proof of Proposition \ref{pr m1} shows that
$$u_k\le w_{t_p}\le (c_6t_pk^p+k) w_0,$$
then
\begin{eqnarray*}
\int_{\Omega} u_{k}^{p-1}\xi^2 \,dx & \le&   \arraycolsep=1pt\left\{
\begin{array}{lll}
 \displaystyle  c_{9} k^{p-1} \int_{\Omega}  \xi^2(x) |x|^{ (p-1)\tau_-(\mu)} \,dx\quad
   &{\rm if}\quad  \mu>\mu_0,\\[1.5mm]
 \phantom{   }
 \displaystyle c_{9} k^{p-1}\int_{\Omega}  \xi^2(x) |x|^{ (p-1)\tau_-(\mu)} (|\ln |x||+1)^{p-1}\,dx   \quad \ &{\rm   if}\quad   \mu=\mu_0,
 \end{array}\right.
\end{eqnarray*}
where $ (p-1)\tau_-(\mu)>-2$ by the assumption $p<p^*_\mu$. By the improved inequality, see \cite[(1.5)]{FT}, we have that
$$\int_\Omega |\nabla \xi|^2 dx \ge \frac{(N-2)^2}4 \int_\Omega |x|^{-2} \xi^2 dx+ c_{10}\int_\Omega |V| \xi^2 dx$$
where   $V\in L^q(\Omega)$ with $q>\frac N2$.
For $V(x)= |x|^{ (p-1)\tau_-(\mu)}$ or $V(x)= |x|^{ (p-1)\tau_-(\mu)} (|\ln |x||+1)^{p-1}$, it is obvious  that
$V\in L^q(\Omega)$ with $q>\frac N2$.
Then if $k>0$  small enough, we have that
\begin{eqnarray*}
p \int_{\Omega} u_{k}^{p-1}\xi^2 \,dx  <    \norm{\xi}^2_\mu,
\end{eqnarray*}
Then $u_{k} $ is a stable solution of (\ref{eq 1.1}) for $k>0$ small.\smallskip

\emph{Proof of the stability for  $k\in(0,k^*)$.}
Suppose that if $u_k$ is not stable for some $k \in (0, k^*)$, then we have that
\begin{equation}\label{e.3.5}
\sigma_1:= \inf_{\xi\in H^1_0(\Omega)\setminus\{0\}}\frac{\norm{\xi}_\mu^2}{p\int_{\Omega}  u_k^{p-1}\xi^2 \,dx}\le 1.
\end{equation}
It is clear that $\sigma_1$ is achieved by some function $\xi_1$, which can be taken as nonnegative and satisfies
$$ \mathcal{L}_\mu   \xi_1 =  \sigma_1 pu_{k}^{p-1}   \xi_1\quad{\rm and}\quad \xi_1\in H^1_0(\Omega)\cap C^2(\Omega\setminus\{0\}).$$
Let $\xi_\epsilon^*$ be a nonnegative function in $C^2_0(\Omega)$  such that
$$\xi_\epsilon^*=\xi_1^*\quad{\rm in}\quad \Omega\setminus B_\epsilon(0)\quad{\rm and}\quad \xi_\epsilon^*\le \xi_1^*\quad{\rm in}\quad   B_\epsilon(0), $$
where $\epsilon>0$ will be determined latter.
Then $\xi_\epsilon=\xi_\epsilon^* \Gamma_\mu$ and
 \begin{equation}\label{eq 3.1}
\mathcal{L}_\mu^*   \xi_\epsilon^* =  \sigma_1 pu_{k}^{p-1}  \xi_\epsilon^*  \qquad{\rm in}\quad\Omega\setminus B_\epsilon(0),
 \end{equation}

Choosing $\hat{k}\in (k,k^*)$ and letting $w=u_{\hat k}- u_{k}$ be the solution of
$$
\arraycolsep=1pt\left\{
\begin{array}{lll}
 \displaystyle  \mathcal{L}_\mu u= u_{\hat k}^p- u_{k}^p\qquad
   &{\rm in}\quad \ {\Omega}\setminus \{0\},\\[1mm]
 \phantom{  L_\mu \, }
 \displaystyle  u= 0\qquad  &{\rm   on}\quad\ \partial{\Omega},\\[1mm]
 \phantom{   }
  \displaystyle \lim_{x\to0}u(x)\Phi_\mu^{-1}(x)=\hat k-k.
 \end{array}\right.
 $$
By the elementary inequality
\begin{equation}\label{e.3.501q}
(a+b)^p\ge a^p+pa^{p-1}b \quad{\rm for}\ \ a, b\ge 0,
\end{equation}
we infers that
\begin{eqnarray*}
\int_\Omega  w \mathcal{L}_\mu^* \xi_{\epsilon}^*\, d\mu-\int_{B_\epsilon(0)}  w \mathcal{L}_\mu^* \xi_{\epsilon}^*\, d\mu &=& \int_\Omega (u_{\hat k}^p- u_{k}^p)\xi_{\epsilon}^*\, d\mu+ (\hat k-k)\xi_{\epsilon}^*(0)-\int_{B_\epsilon(0)}  w \mathcal{L}_\mu^* \xi_{\epsilon}^*\, d\mu \\
    &>& \int_{\Omega}p u_{k}^{p-1} w\xi_{\epsilon}^*\,d\mu -\int_{B_\epsilon(0)}  w \mathcal{L}_\mu^* \xi_{\epsilon}^*\, d\mu.
\end{eqnarray*}

From (\ref{eq 3.1}), we obtain that
\begin{eqnarray*}
\sigma_1\int_{\Omega }pu_{k}^{p-1} w  \xi_{\epsilon}^*\,d\mu&\ge& \sigma_1\int_{\Omega\setminus B_\epsilon(0)}pu_{k}^{p-1} w  \xi_{\epsilon}^*\,d\mu=\int_{\Omega\setminus B_\epsilon(0) }  w \mathcal{L}_\mu^* \xi_\epsilon^* \,d\mu\\&=&\int_\Omega  w \mathcal{L}_\mu^* \xi_{\epsilon}^*\, d\mu-\int_{B_\epsilon(0)}  w \mathcal{L}_\mu^* \xi_{\epsilon}^*\, d\mu.
\end{eqnarray*}
Then
$$\frac{1}{2}(1-\sigma_1) \int_{\Omega} pu_{k}^{p-1} w   \xi_{\epsilon}^* \,d\mu < \int_{B_\epsilon(0)}  w \mathcal{L}_\mu^* \xi_{\epsilon}^*\, d\mu,$$
where the mapping $\epsilon\to\displaystyle\int_{\Omega} pu_{k}^{p-1} w   \xi_{\epsilon}^* \,d\mu $ is decreasing.
Since $$\lim_{x\to0}w(x)\Phi_\mu^{-1}(x)=\hat k-k,$$
then for $\epsilon>0$ small, we have that  $w(x)\le 2(\hat k-k)\Phi_\mu$ in $B_\epsilon(0)\setminus\{0\}$ and
\begin{eqnarray*}
|\int_{B_\epsilon(0)}w \mathcal{L}_\mu^*(\varphi_\epsilon^*)\,d\mu|
 \le   \arraycolsep=1pt\left\{
\begin{array}{lll}
c_{11} \epsilon\quad    &{\rm for}\quad \mu>\mu_0,\\[1mm]
 \phantom{  }
c_{11} \epsilon(-\ln \epsilon)\quad \ \ \quad &{\rm for}\quad \mu=\mu_0,
 \end{array}\right.
\end{eqnarray*}
which implies that
$$\lim_{\epsilon\to0}\int_{B_\epsilon(0)} w \mathcal{L}_\mu^*(\varphi_\epsilon^*)\,d\mu=0,$$
which is impossible. Consequently,
$$
 p\int_{\Omega}    u_{k}^{p-1}  \xi^2 \,dx < \norm{\xi}_\mu^2,\qquad \forall \xi\in H^1_0(\Omega).
$$
and then
$u_{k}$ is stable for $0<k<k^*$.\smallskip

{\it Proof of (\ref{e.3.1}).} For any $k\in (0,k^*),$ let $k'=\frac{k+k^*}{2}>k$ and $l_0=\frac k{k'}<1$, then
we see that the minimal solution $u_{k'}$ of (\ref{eq 1.1}) verifying (\ref{1.2}) with $k'$ being  stable. We observe that
$ l_0u_{k'}\ge l_0^p u_{k'}$,
$$\lim_{x\to0}l_0 u_{k'}(x)\Phi_\mu^{-1}(x)=k $$
and
$$\mathcal{L}_\mu(l_0 u_{k'}) = l_0u_{k'}^p> (l_0u_{k'})^p\quad{\rm in}\quad \Omega\setminus\{0\}. $$
So   $l_0u_{k'}$ is super solution of (\ref{eq 1.1}) verifying (\ref{1.2}) with such $k$,
and $u_k$ is the minimal solution of (\ref{eq 1.1}) verifying (\ref{1.2}) with such $k$, then
$$l_0u_{k'}\ge u_{k},$$
then for $\xi\in  H^1_0(\Omega)\setminus\{0\}$, we have that
\begin{eqnarray*}
0< \norm{\xi}_\mu^2- p\int_{\Omega}   u_{k'}^{p-1}\xi^2 \,dx
& \le&   \norm{\xi}_\mu^2- pl_0^{1-p}\int_{\Omega}   u_{k}^{p-1}\xi^2 \,dx
\\&=&l_0^{1-p}\left[l_0^{p-1}\norm{\xi}_\mu^2- p\int_{\Omega}  u_{k}^{p-1}\xi^2 \,dx\right],
\end{eqnarray*}
thus,
\begin{eqnarray*}
 \norm{\xi}_\mu^2- p\int_{\Omega} u_{k}^{p-1}\xi^2 \,dx
 &=& (1-l_0^{p-1})\norm{\xi}_\mu^2+  \left[l_0^{p-1}\norm{\xi}_\mu^2- p\int_{\Omega} u_{k}^{p-1}\xi^2 \,dx\right]
   \\ & \ge& (1-l_0^{p-1}) \norm{\xi}_\mu^2,
\end{eqnarray*}
which, together with the fact that
 $$ 1-l_0^{p-1}= 1-(\frac{2k}{k+k^*})^{p-1} $$
 implies (\ref{e.3.1}).\hfill$\Box$\smallskip

 Now we would like to approach the weak solution when $k=k^*$ by the minimal solution $u_k$ with $k<k^*$.

\begin{proposition}\label{pr m2}
Assume that $p\in(1,p^*_\mu)$ and $k^*$ is given by Proposition \ref{pr m1},  then
  problem (\ref{eq 1.2}) with $k=k^*$  admits a minimal nonnegative solution $u_{k^*}$, which is semi-stable.

\end{proposition}

\noindent {\bf Proof.}
Let $(\lambda_1,\varphi_1)$ be the first eigenvalue and positive eigenfunction of $\mathcal{L}_\mu$ in $H^1_0(\Omega)$ and  let $\{\varphi_\epsilon^*\}_\epsilon$ be an increasing  sequence of nonnegative functions in $C^2_0(\Omega)$  verifying (\ref{3.4}), then for $k\in(0,k^*)$, we have that from (\ref{3.3})
\begin{eqnarray*}
  \int_\Omega u_k^p\varphi_\epsilon^*\,d\mu=  (1+\epsilon) \int_\Omega u_k\mathcal{L}_\mu^*\varphi_\epsilon^*\, d\mu- k\varphi_\epsilon^*(0)<     \lambda_1 \left(\int_\Omega u_k^p\varphi_\epsilon^*\,d\mu\right)^{\frac1p} \left(\int_\Omega  \varphi_\epsilon^*\,d\mu\right)^{1-\frac1p},
\end{eqnarray*}
which implies that
\begin{equation}\label{3.2}
 \norm{u_k}_{L^p(\Omega, \,d\mu)}\le \norm{ \varphi_\epsilon^*}_{L^\infty(\Omega)}^{-1} \lambda_1^{\frac{p}{p-1}}\int_\Omega  \varphi_\epsilon^*\,d\mu.
\end{equation}
Since the mapping $k\mapsto u_k$ is increasing, then  the limit of $\{u_k\}_k$ exists, denoting $u_{k^*}$,  $u_k\to u_{k^*}$ in $  L^p(\Omega, \,d\mu)$ and
$$\int_{\Omega} u_{k^*}\mathcal{L}_\mu^*  \xi \,d\mu =\int_{\Omega} u_{k^*}^p\xi \,d\mu +{k^*}\xi(0), \quad \forall\, \xi\in C^2_c(\Omega).$$
So we conclude that (\ref{eq 1.2}) has a weak solution and then (\ref{eq 1.2}) has minimal solution $u_{k^*}$.

{\it Proof of the semi-stability of $u_{k^*}$.}
For any $\epsilon>0$ and $\xi\in H^1_0(\Omega)\setminus\{0\}$, there exists $k(\epsilon)>0$ such that for all $k\in(k(\epsilon),k^*)$,
\begin{eqnarray*}
p\int_{\Omega} u_{k^*}^p\xi \,dx   \le   p\int_{\Omega} u_{k}^p\xi \,dx    +
  (k^*-k) p\int_{\Omega} u_{k^*}^p\xi \,dx \le   \norm{\xi}^2_{\mu} +\epsilon.
\end{eqnarray*}
By the arbitrary of $\epsilon>0$, we have that
$u_{k^*}$ is semi-stable. \hfill$\Box$

\vskip0.3cm
In the special case that $\Omega=B_1(0)$, we have a complete analysis on the minimal solution of problem (\ref{eq 1.1}) subjecting to (\ref{1.2}).

\begin{corollary}\label{cr m1}
Assume that $\Omega=B_1(0)$ and $p\in(1,p^*_\mu)$,  then there exists $k^*>0$
  such that

 $(i)$ for $ k\in(0, k^*)$, problem (\ref{eq 1.1}) subjecting to (\ref{1.2})  admits a minimal nonnegative solution $u_{k}$, which is stable.  Furthermore, $u_k$ is a $d\mu$-distributional solution of (\ref{eq 1.2}).

$(ii)$
for $k> k^*$,  problem (\ref{eq 1.1}) subjecting to (\ref{1.2}) or (\ref{eq 1.2}) admits no solution.

 $(iii)$ for $ k=k^*$, problem (\ref{eq 1.1}) subjecting to (\ref{1.2}) with $k=k^*$  admits a  unique  nonnegative solution $u_{k^*}$, which is semi-stable, and $u_{k^*}$ is a $d\mu$-distributional solution of (\ref{eq 1.2}).

\end{corollary}
{\bf Proof.} From Proposition \ref{pr m1}, Lemma \ref{lm stability} and Proposition \ref{pr m2}, we only have to prove that problem (\ref{eq 1.1}) subjecting to (\ref{1.2}) with $k=k^*$
has a unique solution $u_{k^*}$, which is semi-stable.

The minimal solution $u_k$ of problem (\ref{eq 1.1}) subjecting to (\ref{1.2}) is approximated by $v_n$, the solution of (\ref{2.02-2}). It is obvious that $v_n$ is radially symmetric and non-increasing  with respect to $|x|$, so is $u_k$.

We observe that the mapping $k\mapsto u_k$ is increasing and $u_{k^*}:=\lim_{k\to k^*} u_k$ a.e. in $B_1(0)$ and $u_{k^*}\in L^p(B_1(0),d\mu)$, thus, $u_{k^*}$ is radially symmetric and non-increasing  with respect to $|x|$. Now we claim that $u_{k^*}$ is locally bounded in $B_1(0)\setminus\{0\}$. In fact, if there exists $x_0\not=0$ such that $u_{k^*}(x_0)=+\infty$, then $u_{k^*}=+\infty$ in $\overline{B_{|x_0|}(0)}$, which is impossible with $u_{k^*}\in L^p(B_1(0),d\mu)$.

Denote
$$k^{**}=\sup\{k>0:\ (\ref{eq 1.2}) \ {\rm  has\ minimal\ solution} \}.$$
Obviously, we have that $k^{**}\ge k^*$. We show that
\begin{equation}\label{k**}
 k^{**}=k^*.
\end{equation}
 In fact, for any $ k< k^{**}$, let $w_k$ be a nonnegative solution of (\ref{eq 1.2}),
then  problem (\ref{2.02-2}) has solution  $v_n$, which is a $d\mu$ distributional solution of
$$\mathcal{L}_\mu v_n= v_{n-1}^p+k\delta_0\quad {\rm in}\quad B_1(0),\qquad v_n=0\quad {\rm on}\quad \partial B_1(0), $$
where $v_0=kV_0$ is the solution of
\begin{equation}\label{eq 1.1k}
 \mathcal{L}_\mu V_0= \delta_0\quad {\rm in}\quad B_1(0),\qquad V_0=0\quad {\rm on}\quad \partial B_1(0).
\end{equation}
Note that $w_k$ is an upper bound of $\{v_n\}$, so the limit of $\{v_n\}$ is the minimal solution of (\ref{eq 1.2}).
From the symmetry and monotonicity of $v_n$ implies that $w_k$ is a classical solution of problem (\ref{eq 1.1}) subjecting to (\ref{1.2}).
So we have that $k^*\ge k^{**}$. As a consequence, (\ref{k**}) holds true.

From   the monotonicity of the mapping $k\mapsto u_k$,  the convergence $u_{k^*}=\lim_{k\to k^*} u_k$ and uniformly bounded locally  in $B_1(0)\setminus\{0\}$, so
it implies by the inner regularity results, we have that  $\{u_k\}$ converges to $ u_{k^*}$ in $C^2(B_1(0)\setminus\{0\})$. 
Furthermore, $u_{k^*}$ is a very weak solution of (\ref{eq 1.2}) with $k=k^*$, which implies that $u_{k^*}$ verifies (\ref{1.2}) with $k=k^*$.

{\it We prove the uniqueness for $k=k^*$.}
 Since $u_{k^*}$ is semi-stable, we have that
$$
\sigma_1:= \inf_{\xi\in H^1_0(\Omega)\setminus\{0\}}\frac{\norm{\xi}^2_{\mu}}{p\int_{\Omega}  u_{k^*}^{p-1}\xi^2 \,dx}\ge 1.
$$

We prove $\sigma_1=1$. If not, we may assume that $\sigma_1>1$. We note that the minimal solution $u_k$ could be written as
$$u_k=kV_0+w_k,$$
where $V_0$ is the solution of (\ref{eq 1.1k}) and $w_k\in H^1_0(B_1(0))$ is a solution of
$$\mathcal{L}_\mu w_k =  (kV_0+w_k)^p\quad {\rm in}\quad B_1(0). $$

 Denote
$E:(0,+\infty)\times H^1_0(B_1(0))\to H^{-1}(B_1(0))$ as
$$E(k,u)=\mathcal{L}_\mu u-  (kV_0+u)^p\quad {\rm and}\quad E_{u}(k,w_k)w=\mathcal{L}_\mu w-  pu_k^{p-1} w.$$
Since
\begin{eqnarray*}
\langle E_u(k^*,w_{k^*})w,\, w \rangle  &:=&\int_{B_1(0)}\left( |\nabla w|^2+\frac{\mu}{|x|^2}w^2-  pu_{k^*}^{p-1} w^2\right) dx \\
    &\ge & (1-\frac1{\sigma_1}) \norm{w}_{\mu}^2.
\end{eqnarray*}
So for any $f\in H^{-1}(B_1(0))$, there exists a unique solution $w_f$ in $H^1_0(B_1(0))$ such that
$$E_u(k^*,w_{k^*}) w_f=f,$$
that is, $E_u(k^*,w_{k^*}): H^1_0(B_1(0))\to H^{-1}(B_1(0))$ is invertible.  Then, by the implicit function theorem, there exists $\epsilon > 0$ such that
$E(k,w_k) = 0$ has a solution $w_k\in H^1_0(B_1(0))$ for $k\in (k^*-\epsilon,\,k^*+\epsilon)$, then $kV_0+w_k$ is a $d\mu$-distributional solution of (\ref{eq 1.2}) with
$k>k^{**}$ by the fact $k^*=k^{**}$. This contradicts the definition of $k^{**}$. Thus,  $\sigma_1=1$.

 By the compact embedding theorem, $\sigma_1=1$ is achievable and  its achieved function $\xi_1$ could be setting to be nonnegative and   satisfies
$$ \mathcal{L}_\mu   \xi_1   =   pu_{k^*}^{p-1}   \xi_1\ \ {\rm in}\ \ B_1(0), \qquad \xi=0\ \ {\rm on}\ \ \partial B_1(0).$$

If problem (\ref{eq 1.1}) admits a solution $u>u_{k^*}$.
Let $w=u-u_{k^* }>0$, then it verifies
$$
\arraycolsep=1pt\left\{
\begin{array}{lll}
 \displaystyle  \mathcal{L}_\mu w= u^p- u_{k^*}^p\qquad
   &{\rm in}\quad \ {B_1(0)}\setminus \{0\},\\[1mm]
 \phantom{  L_\mu \, }
 \displaystyle  w= 0\qquad  &{\rm   on}\quad\ \partial{B_1(0)},\\[1mm]
 \phantom{   }
  \displaystyle \lim_{x\to0}w(x)\Phi_\mu^{-1}(x)=0.
 \end{array}\right.
 $$
Then
\begin{eqnarray*}
 \int_{B_1(0)}pu_{k^*}^{p-1} w\xi_1\,dx   =     \int_{B_1(0)} w \mathcal{L}_\mu^*\xi_1\,dx=\int_{B_1(0)}  \xi_1(u^p-u_{k^*}^p)\,dx
 >   \int_{B_1(0)}p u_{k^*}^{p-1} w\xi_1\,dx,
\end{eqnarray*}
which is impossible.  As a conclusion, $u_{k^*}$ is the unique solution of (\ref{eq 1.1}) with $k=k^*$.\hfill$\Box$

\subsection{Mountain Pass Solution}

For the second solution of (\ref{eq 1.1}), we would like to apply the Mountain-Pass theorem to find a positive weak solution of
\begin{equation}\label{eq 4.1}
  \arraycolsep=1pt\left\{
\begin{array}{lll}
 \displaystyle   \mathcal{L}_\mu  u= (u_{k}+u_+)^{p}- u_{k}^p \quad
 &{\rm in}\quad \Omega,\\[2mm]
 \phantom{  \mathcal{L}_\mu }
 u \in H^1_0(\Omega),
\end{array}\right.
\end{equation}
where $k\in (0,k^*)$ and $u_{k}$ is the minimal positive solution of (\ref{eq 1.1}) obtained by Thoerem \ref{teo 0}.
The second solution of (\ref{eq 1.1}) is derived by following proposition.

\begin{proposition}\label{pr 3.1}
Assume that   $p\in(1,p_\mu^*)$, $k\in (0,k^*)$   and
$u_{k}$ is the minimal positive solution of (\ref{eq 1.1}) subject to (\ref{1.2}).
Then problem (\ref{eq 4.1}) has a positive solution $v_k$ satisfying that $v_k>u_{k}$ in $\Omega$ and
\begin{equation}\label{eq 4.2}
  \arraycolsep=1pt\left\{
\begin{array}{lll}
 \displaystyle   \mathcal{L}_\mu  u= (u_{k}+u_+)^{p}- u_{k}^p \quad
 &{\rm in}\quad \Omega\setminus\{0\},\\[2mm]
 \phantom{  \mathcal{L}_\mu }
 u=0  &{\rm on}\quad \partial\Omega,\\[2mm]
 \phantom{   }
 \displaystyle \lim_{x\to0}u(x)\Phi_\mu^{-1}(x)=0.
\end{array}\right.
\end{equation}

\end{proposition}
 {\bf Proof.}
We would like to employe the Mountain Pass Theorem to look for the weak solution of (\ref{eq 4.1}).  A function $v$ is said to be a  weak solution of (\ref{eq 4.1}) if
\begin{equation}\label{weak 1}
\langle u,\xi\rangle_\mu  =\int_{\Omega}\left[ (u_{k}+u_+)^{p}- u_{k}^p\right]\xi \,dx,\quad \forall\xi\in  H^1_0(\Omega).
\end{equation}
The natural functional associated to (\ref{eq 4.1}) is the following
\begin{equation}\label{E}
E(v)=\frac12\norm{v}_\mu^2-\int_{\Omega}F(u_{k},v_+)\,dx, \quad\forall v\in  H^1_0(\Omega),
\end{equation}
where
\begin{equation}\label{F}
 F(s,t)=\frac1{p+1}\left[(s+t_+)^{p+1}-s^{p+1}-(p+1)s^pt_+\right].
\end{equation}

We observe that for any $\epsilon>0$, there exists some $c_{\epsilon}>0$, depend only on $p$, such that
$$0\le  F(s,t) \le (p+\epsilon)s^{p-1}t^2+c_{\epsilon}t^{p+1},\quad s,t\ge0,$$
then for any $v\in  H^1_0(\Omega)$, we have that
\begin{eqnarray*}
 \int_{\Omega}F(u_{k},v_+) \,dx  &\le& (p+\epsilon)\int_{\Omega}u_{k}^{p-1}v_+^2 \,dx +c_{\epsilon}\int_{\Omega}v_+^{p+1} \,dx  \\
    &\le& c_{12}\norm{v}_\alpha^2,
\end{eqnarray*}
thus, $E$ is well defined in $  H^1_0(\Omega)$.

We observe that $E(0)=0$ and let $v\in  H^1_0(\Omega)$ with $\norm{v}_{ \alpha}=1$, then
 for  $k\in(0,k^*)$, choosing $\epsilon>0$ small enough, it infers from (\ref{e.3.1}) that
\begin{eqnarray*}
   E(tv) &= & \frac12t^2\norm{v}_\mu^2- \int_{\Omega}F(u_{k},tv_+)\,\,dx \\
   &\ge&t^2\left(\frac12\|v\|_\mu^2- (p+\epsilon) \int_{\Omega}  v_k^{p-1} v^2 \,dx\right)- c_{13} t^{p+1}\int_{\Omega}   |v|^{p+1} \,dx
   \\
   &\ge &c_{14} t^2\|v\|_\mu^2 - c_{15} t^{p+1} \|v\|_\mu^{p+1}
    \\
   &\ge & \frac{1}2c_{14}t^2-c_{15}t^{p+1},
\end{eqnarray*}
where  we used   (\ref{3.1}) in the first inequality.
Then there exists   $\sigma_0>0$  small such that  for $\|v\|_{ H^1_0(\Omega)} =1$,  we have that
$$ E(\sigma_0 v)\ge    \frac{c_{14}} 4\sigma_0^2=:\beta>0.$$

 We take a nonnegative function $v_0\in  H^1_0(\Omega)$,  then
 $$ F(u_{k},tv_0)\ge \frac{1}{p+1}t^{p+1} v_0^{p+1}-tu_k^pv_0.$$
  Since the space of $\{tv_0: t\in\R\}$ is
a subspace of $ H^1_0(\Omega)$ with dimension 1 and all the norms are equivalent, then $\displaystyle\int_{\Omega} v_0(x)^{p+1} \,dx>0$.
Then there exists $t_0>0$ such that for $t\ge t_0$,
\begin{eqnarray*}
  E(tv_0)&=& \frac {t^2}2\|v_0\|_\mu^2- \int_{\Omega}   F(u_{k},tv_0)  \,dx\\
   &\le &  \frac {t^2}2\|v_0\|_\mu^2-c_{16}t^{p+1} \int_{\Omega}   v_0^{p+1} \,dx+t\int_{\Omega}   u_{k}^{p}v_0 \,dx
  \\ &\le &c_{17}(t^2+t-t^{p+1})\le  0.
\end{eqnarray*}
 Choosing $e=t_0v_0$, we have that $E(e)\le0$.

We next prove that $E$ satisfies $(PS)_c$ condition. We say that $ E$ has $(PS)_c$ condition, if for any sequence
$\{v_n\}$ in $ H^1_0(\Omega)$ satisfying $ E(v_n)\to c$ and $ E'(v_n)\to0$ as $n\to\infty$, there
is a convergent subsequence. Here the energy level $c$ of functional $E$ is characterized by
\begin{equation}\label{c}
 c=\inf_{\gamma\in \Upsilon}\max_{s\in[0,1]}E(\gamma(s)),
\end{equation}
where $\Upsilon =\{\gamma\in C([0,1]:\, H^1_0(\Omega)):\gamma(0)=0,\ \gamma(1)=e\}$. We observe that
$$c\ge \beta.$$

Let $\{v_n\}$ in $ H^1_0(\Omega)$ satisfying $ E(v_n)\to c$ and $ E'(v_n)\to0$ as $n\to\infty$, then we only have to show that there are a subsequence, still denote it by  $\{v_{n}\}$
 and $v\in  H^1_0(\Omega)$ such that
$$v_{n}\to v\quad {\rm in}\ \ L^2(\Omega,  u_{k}^{p-1}\,dx)\quad{\rm and}\ \ L^{p+1}(\Omega) \quad {\rm as}\ n\to\infty.$$
 We see that
$$
 c_{18}\|w\|_\mu\ge E'(v_n)w
   =   \langle v_n, w\rangle_{\mu} -\int_{\Omega}  f(u_{k}, (v_n)_+)w\,dx
$$
and
\begin{eqnarray}\label{4.2}
 c+1\ge E(v_n) = \frac12\|v_n\|^2_\mu - \int_{\Omega}  F(u_{k}, (v_n)_+)\,dx.
\end{eqnarray}

Let $c_p=\min\{1,p-1\}$, then it follows by \cite[C.2 $(iv)$]{NS}  that
$$f(s,t)t-(2+c_p)F(s,t)\ge -\frac{c_pp}{2}s^{p-1}t^2,\quad s,t\ge0,$$
thus, we multiply (\ref{4.2}) by $(2+c_p)$ and minus $\langle E'(v_n),(v_n)_+\rangle$, we get that
\begin{eqnarray*}
c+ c_{18}\|v_n\|_\mu&\ge&  \frac{c_p}2\|v_n\|^2_\mu -\int_{\Omega} \left[(2+c_p)F(u_{k},(v_n)_+)-f(u_{k},(v_n)_+)(v_n)_+ \right]\,dx   \\
    &\ge &  \frac{c_p}2 \left[\|v_n\|^2_\mu-p \int_{\Omega}  u_{k}^{p-1}v_n^2\,dx\right]
    \\&\ge&  c_{19} \frac{c_p}2\|v_n\|^2_\mu.
\end{eqnarray*}
Then we derive that $v_n$ is uniformly bounded in $ H^1_0(\Omega)$ for $k\in(0,k^*)$.
Thus,  there exists a subsequence $\{v_{n}\}$ and $v$ such that
$$v_{n}\rightharpoonup v\quad{\rm in}\quad  H^1_0(\Omega), $$
$$v_{n}\to v\quad{\rm a.e.\ in}\ \Omega\quad{\rm and\ \ in}\quad L^{p+1}(\Omega),\ \ L^{2}(\Omega,  u_{k}^{p-1}\,dx),   $$
when $n\to \infty$. Here we have used that
\begin{equation}\label{4.1}
 u_{k}(x)^{p-1}\le c_{20}(1+|\Phi_\mu|)^{p-1} (x),\quad x\in\Omega\setminus\{0\},
\end{equation}
where $|\Phi_\mu (x)|^{p-1}=|x|^{(p-1)\tau_-(\mu)}$ if $\mu>\mu_0$ or $|\Phi_\mu (x)|^{p-1}=|x|^{(p-1)\tau_-(\mu)}|\ln|x||^{p-1}$ if $\mu=\mu_0$,  $( (p-1)\tau_-(\mu)>-2$ and from \cite{NS}, the embedding: $H^1_0(\Omega)\hookrightarrow L^q(\Omega, (1+|\Phi_\mu|)^{p-1}\,dx)$
is compact for $q\in[1, \frac{2N+2(-\tau_-(\mu))(p-1)}{N-2 })$,  particularly, for $q=2$.

We observe that
\begin{eqnarray*}
&&|F(u_{k},v_n)- F(u_{k},v)|
\\&&\qquad= \frac1{p+1}|(u_{k}+(v_n)_+)^p-(u_{k}+v_+)^p-(p+1)u_{k}^p((v_n)_+-v_+)| \\
   &&\qquad\le c_{21}u_{k}^{p-1} ((v_n)_+-v_+)^2+c_{31}((v_n)_+-v_+)^{p+1},
\end{eqnarray*}
which implies that
$$F(u_{k},v_n)\to F(u_{k},v)\quad{\rm a.e.\ in}\ \Omega\quad{\rm and\ \ in}\quad L^1(\Omega).$$
Then, together with $\lim_{n\to\infty} E(v_{n})=c$,
we have that $\|v_{n}\|_\mu\to \|v\|_\mu$  and
$v_{n}\to v$ in $ H^1_0(\Omega)$ as $n\to\infty$.

Now Mountain Pass Theorem  (for instance, \cite[Theorem 6.1]{struwe}; see also \cite{rabinowitz}) is applied to obtain that
there exists a critical point $v\in  H^1_0(\Omega)$
 of $ E$ at some value $c\ge \beta>0$. By $\beta>0$, we have that $v$ is nontrivial and nonnegative. Then $v$ is a positive weak solution
of (\ref{eq 4.1}).  By using bootstrap argument in \cite{HL}, the interior regularity of $v$ could be improved to be in $ H^1_0(\Omega)\cap C^2(\Omega\setminus\{0\})$, since $u_{k}$ is locally bounded in $\Omega\setminus\{0\}$ and $p<p_\mu^*$. From (\ref{4.1}) with $p<p_\mu^*$, we have that there is some $q>\frac N{2}$ such that
$$ u_{k}^{p-1}\in L^q(\Omega),$$
then $v_k$ is bounded at the origin. Therefore, we have that $v_k$ is a solution of (\ref{eq 4.2}).  Moreover, by Maximum Principle, we conclude that $v>0$ in  $\Omega$.
 \hfill$\Box$

\medskip
\noindent{\it Proof of  Theorem \ref{teo 1}.} Proposition \ref{pr m1} shows the existence and nonexistence of  minimal singular solution of
(\ref{eq 1.1}) subjecting to (\ref{1.2}) with $k\in(0,k^*)$.  From Proposition \ref{pr 3.1}, we obtain that there is a  positive weak solution of $v_k$
of (\ref{eq 4.2}),
then $(u_{k}+v_k)$ satisfies
$$
\mathcal{L}_\mu (u_{k}+v_k)  =   (u_{k}+v_k)^{p}\quad {\rm and}\quad \lim_{x\to0}(u_{k}+v_k)(x)\Phi_\mu^{-1}(x)=0,
$$
which means that $v_k +u_{k}$ is a classical solution of  (\ref{eq 1.1}) subjecting to (\ref{1.2}).

Therefore, Theorem \ref{teo 1} part $(i)$ and part $(iii)$ hold.

For the extremal case, Theorem \ref{teo 1} part $(ii)$ follows by Proposition   \ref{pr m2}  and Corollary \ref{cr m1} in the case of $\Omega=B_1(0)$.
\hfill$\Box$

\bigskip\bigskip

 \noindent{\bf Acknowledgements:} H. Chen  is supported by NNSF of China, No:11401270,   by SRF for ROCS, SEM and by the Jiangxi Provincial Natural Science Foundation, No: 20161ACB20007.
   F. Zhou is partially supported by NNSF of China, No:11271133 and No:11431005,
 and Shanghai Key Laboratory of PMMP.

\end{document}